\def\withaddnotes{0}
\def\withffootnotes{0}
\def\withrefnotes{0}

\newcommand{\addnote}[1]{\ifnum\withaddnotes=0\else{\footnote{#1}}\fi}
\newcommand{\ffootnote}[1]{\ifnum\withffootnotes=0\else{\footnote{#1}}\fi}
\newcommand{\refnote}[1]{\ifnum\withrefnotes=0\else{\footnote{#1}}\fi}

\documentclass[11pt]{article}

\usepackage[T1]{fontenc}
\usepackage[utf8]{inputenc}
\usepackage{amsmath, amssymb, amsfonts}
\usepackage{array}
\usepackage{float}
\usepackage{cancel}
\usepackage{graphicx}
\usepackage{amsthm}
\usepackage{stmaryrd}
\usepackage[colorlinks,linkcolor=blue,citecolor=blue,pagebackref]{hyperref}

\def\mat#1#2{%
  \if#1<%
    \if#2=\preccurlyeq\else\prec#2\fi%
  \else\if#1>%
         \if#2=\succcurlyeq\else\succ#2\fi%
       \else#1#2%
       \fi%
  \fi}

\renewcommand{\H}{^\mathsf{H}}
\newcommand{\Kb}{\textbf{K}}
\newcommand{\NN}{\mathbb{N}}
\newcommand{\NNn}{\NN^n}

\newcommand{\NNp}{\NN^p}
\newcommand{\RR}{\mathbb{R}}
\newcommand{\RRn}{\RR^n}
\newcommand{\RRnn}{\RR^{2n}}

\newcommand{\calA}{\mathcal{A}}

\newcommand{\eb}{\textbf{e}}

\newcommand{\tr}{\operatorname{tr}}
\newcommand{\vb}{\textbf{v}}
\newcommand{\xb}{\textbf{x}}
\newcommand{\yb}{\textbf{y}}

\renewcommand{\sb}{\textbf{s}}

\newcommand{\T}{^\mathsf{T}}

\newcommand{\dmu}{{\rm d}\mu}

\let\leq\leqslant
\let\geq\geqslant

\newcommand{\thetitle}{%
Application of the Moment--SOS Approach to Global Optimization of the
OPF Problem\ifnum\withrefnotes=0\else\\(referee version)\fi
}

\newcommand{\theabstract}{%
Finding a global solution to the optimal power flow (OPF) problem is
difficult due to its nonconvexity. A convex relaxation in the form of
semidefinite programming (SDP) has attracted much attention lately as it
yields a global solution in several practical cases. However, it does
not in all cases, and such cases have been documented in recent
publications. This paper presents another SDP method known as the
moment-sos (sum of squares) approach, which generates a sequence that
converges towards a global solution to the OPF problem at the cost of
higher runtime. Our finding is that in the small examples where the
previously studied SDP method fails, this approach finds the global
solution. The higher cost in runtime is due to an increase in the matrix
size of the SDP problem, which can vary from one instance to another.
Numerical experiment shows that the size is very often a quadratic
function of the number of buses in the network, whereas it is a linear
function of the number of buses in the case of the previously studied
SDP method.
}

\newcommand{\thekeywords}{%
Global optimization,
moment/sum-of-squares approach,
optimal power flow,
polynomial optimization,
semidefinite programming.
}
\begin{document}

\pagestyle{plain}

\renewcommand{\thefootnote}{\fnsymbol{footnote}}
\setcounter{footnote}{0}

\begin{center}
{\Large\bf\thetitle\footnote{This work was partly supported by the
contract CIFRE number 2013/0179 of the ANRT (Association Nationale de la
Recherche et de la Technologie, France).}\\}
\vspace{2ex}\rm\normalsize
C. {\sc Josz}\,%
\footnote{French transmission system operator \textit{Réseau de
Transport d'Electricité} (RTE), 9, rue de la Porte de Buc, BP 561,
F-78000 Versailles, France. E-mails:
\texttt{firstname.lastname@rte-france.com}.}%
\footnote{INRIA Paris-Rocquencourt, BP~105, F-78153~Le Chesnay,
France. E-mails: \texttt{Cedric.\linebreak[2]Josz@\linebreak[2]inria.\linebreak[2]fr},
\texttt{Jean-Charles.\linebreak[2]Gilbert@\linebreak[2]inria.fr}.},
\addtocounter{footnote}{-2}
J.\ {\sc Maeght}\,\footnotemark,
\addtocounter{footnote}{-1}
P.\ {\sc Panciatici}\,\footnotemark,
and
J.\ Ch.\ {\sc Gilbert}\,\footnotemark

\vspace*{0.5ex}
November 18, 2013

\end{center}

\renewcommand{\thefootnote}{\arabic{footnote}}
\setcounter{footnote}{0}

\vspace{-1ex}
\begin{quote}
\small
\theabstract

\vspace{2ex}
{\bf Keywords:}
\thekeywords
\end{quote}

\section{Introduction}
\indent The optimal power flow (OPF) gives its name to a problem pertaining to power systems that was first introduced by Carpentier in 1962
\cite{carpentier.m.j-1962}. It seeks to determine a steady state
operating point of an alternating current (AC) power network that is
optimal under some criteria such as generating costs. The problem can be
cast as a nonlinear optimization problem, which is NP-hard, as was shown
in \cite{lavaei-low-2012}. So far, the various methods
\cite{huneault-galiana-1991, pandya-joshi-2008} that have been
investigated to solve the OPF can only guarantee local optimality, due
to the nonconvexity of the problem. Recent progress suggests that it may
be possible to design a method, based on semidefinite programming (SDP),
that yields global optimality rapidly.

SDP is a subfield of convex conic
optimization~\cite{saigal-vandenberghe-wolkowicz-2000}. It deals with problems whose structure resembles
that of a linear optimization problem, but where the variable that is being
solved for is a positive semidefinite matrix.
An SDP problem has a convex feasible set whose definition is sufficiently general to
model a large variety of convex problems. Furthermore, it can be solved
by efficient techniques, notably the interior point methods, which are
able to find a solution of a given precision in polynomial time. These properties
make the SDP modelling adapted to many applications~\cite{bental-nemirovski-2001}.

The first attempt to use SDP to solve the OPF problem was made by Bai et
al. \cite{bai-fujisawa-wang-wei-2008} in 2008. In
\cite{lavaei-low-2012}, Lavaei and Low show that the OPF can be written
as an SDP problem, with an additional constraint imposing that the rank
of the matrix variable must not exceed~1.
They discard the rank constraint, as it is done in Shor's relaxation
\cite{shor-1987b}, a well-known procedure which applies to quadratically
constrained quadratic problems (see~\cite{sturm-zhang.s-2003,
luo.z.q-ma.w.k-so.a-ye.y-zhang.s-2010} and the references therein).
They also accept quartic terms that appear in some formulations of the
OPF, transforming them by Schur's complement.
Their finding is that
for all IEEE benchmark networks, namely the 9, 14, 30, 57, 118, and
300-bus systems, the rank constraint is satisfied if a small resistance
is added in the lines of the network that have zero resistance. Such a
modification to the network is acceptable because in reality, resistance
is never equal to zero.

There are cases when the rank constraint is not satisfied and a global
solution can thus not be found. Lesieutre et al.
\cite{borden-demarco-lesieutre-molzahn-2011} illustrate this with a
practical 3-bus cyclic network. Gopalakrishnan et al.
\cite{biegler-gopalakrishnan-nikovski-raghunathan-2012} find yet more
examples by modifying the IEEE benchmark networks. Bukhsh et al.
\cite{bukhsh-grothey-mckinnon-trodden-2013} provide a 2-bus and a 5-bus
example. In addition, they document the local solutions to the OPF in
many of the above-mentioned examples where the rank constraint is not
satisfied \cite{bukhsh-grothey-mckinnon-trodden-2013b}.

Several papers propose ways of handling cases when the rank constraint
is not satisfied. Gopalakrishnan et al.
\cite{biegler-gopalakrishnan-nikovski-raghunathan-2012} propose a branch
and reduce algorithm. It is based on the fact that the rank relaxation
gives a lower bound of the optimal value of the OPF. But according to
the authors, using the classical Lagrangian
dual to evaluate a lower bound is about as
efficient. Sojoudi and Lavaei \cite{lavaei-sojoudi-2012} prove that if
one could add controllable phase-shifting transformers to every loop in
the network and if the objective is an increasing function of generated
active power, then the rank constraint is satisfied. Though numerical
experiments confirm this \cite{farivar-low-2013}, such a modification to
the network is not realistic, as opposed to the one mentioned earlier.

Cases where the rank constraint holds have been identified. Authors of
\cite{bose-chandy-gayme-low-2011, tse-zhang-2011, lavaei-sojoudi-2012b}
prove that the rank constraint is satisfied if the graph of the network
is acyclic and if load over-satisfaction is allowed. This is typical of
distribution networks but it is not true of transmission networks.

This paper examines the applicability of the moment-sos (sum of
squares) approach to the OPF. This approach~\cite{lasserre-2000,
parrilo-2000b, lasserre-2001} aims at finding global solutions to
polynomial optimization problems, of which the OPF is a particular
instance. The approach can be viewed as an extension of the SDP method
of \cite{lavaei-low-2012}. Indeed, it proposes a sequence of SDP
relaxations whose first element is the rank relaxation in many cases.
The subsequent relaxations of the sequence become more and more
accurate. When the rank relaxation fails, it is therefore natural to see
whether the second order relaxation provides the global minimum, then
the third, and so on.

The limit to this approach is that the complexity of the relaxations
rapidly increases. The matrix size of the SDP relaxation of order $d$ is
roughly equal to the number of buses in the network to the power $d$.
Surprisingly, in the 2, 3, and 5-bus systems found in
\cite{borden-demarco-lesieutre-molzahn-2011,
bukhsh-grothey-mckinnon-trodden-2013} where the rank relaxation fails,
the second order relaxation nearly always finds the global solution. 

This paper is organized as follows. Section \ref{sec:Link between
optimal power flow (OPF) and polynomial optimization} presents a
formulation of the OPF problem and shows that it can be viewed as a
polynomial optimization problem. The moment-sos approach which aims at solving such problems is described in section
\ref{sec:moment-sos-approach}. In section \ref{sec:Numerical results},
numerical results show that this approach successfully
finds the global solution to the 2, 3, and 5-bus systems mentioned
earlier. Conclusions are given in section \ref{sec:Conclusion}.

\section{OPF as a polynomial optimization problem}
\label{sec:Link between optimal power flow (OPF) and polynomial optimization}

We first present a classical formulation of the OPF with quadratic
objective, Kirchoff's laws, Ohm's law, power balance equations, and
operational constraints. It allows for ideal phase-shifting transformers
that have a fixed ratio. Next we show how the OPF can be cast as a
polynomial optimization problem.

\subsection{Classical formulation of the OPF}
\label{subsec:Classical formulation of the OPF}

Let $\text{j}$ denote the imaginary unit and let $|z|$ and $z\H$
respectively denote the modulus and the conjugate of a complex number
$z$.

Consider an AC electricity transmission network defined by a set of
buses $\mathcal{N} = \{1,\hdots,n\}$ of which a subset $\mathcal{G}
\subset \mathcal{N}$ is connected to generators. Let $s_k^\text{gen} =
p_k^\text{gen} + \text{j} q_k^\text{gen} \in \mathbb{C}$ denote
generated power at bus $k \in \mathcal{G}$. All buses are connected to a
load (i.e., power demand). Let $s_k^\text{dem} = p_k^\text{dem} + \text{j}
q_k^\text{dem} \in \mathbb{C}$ denote power demand at bus $k \in
\mathcal{N}$. Let $v_k \in \mathbb{C}$ denote voltage at bus $k\in
\mathcal{N}$ and $i_k \in \mathbb{C}$ denote current injected into the
network at bus $k\in \mathcal{N}$. The convention used for current means
that $v_k i_k\H$ is the power injected into the network at bus $k\in
\mathcal{N}$. This means that $v_k i_k\H = -s_k^\text{dem}$ at bus $k
\in \mathcal{N} \setminus \mathcal{G}$ and $v_k i_k\H = s_k^\text{gen} -
s_k^\text{dem}$ at bus $k\in \mathcal{G}$.

\begin{figure}[H]
  \centering
    \includegraphics[width=.5\textwidth]{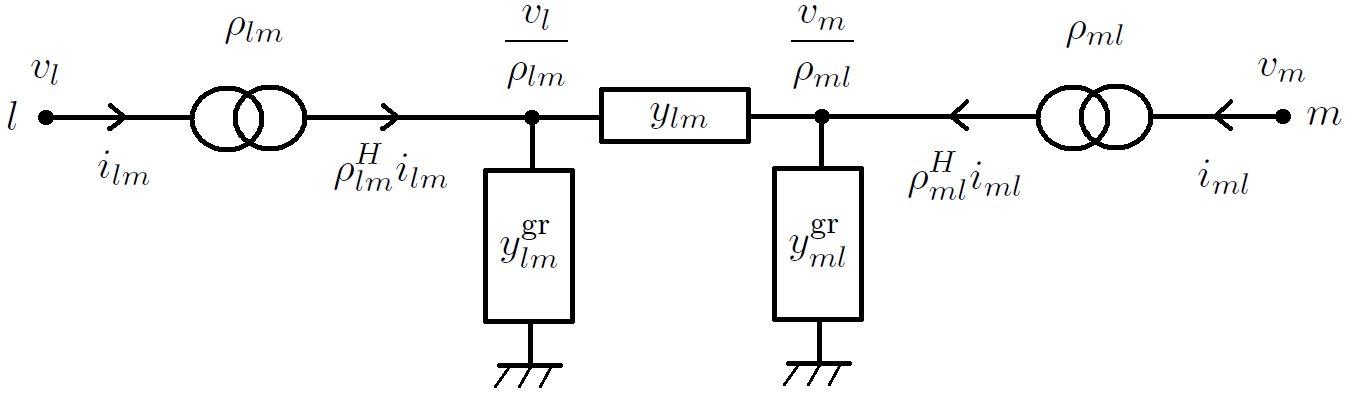}
  \caption{Branch connecting buses $l$ and $m$}
  \label{fig:transmission line with pst}
\end{figure}

The network connects buses to one another through a set of branches
$\mathcal{L} \subset \mathcal{N} \times \mathcal{N}$. Let
$\mathcal{N}(l)$ denote the set of buses connected to bus $l \in
\mathcal{N}$ by a branch in $\mathcal{L}$. If there is a branch
connecting buses $l\in \mathcal{N}$ and $m\in\mathcal{N}$, then $(l,m)
\in \mathcal{L}$ and $(m,l) \in \mathcal{L}$. A branch between two buses
is described in figure \ref{fig:transmission line with pst}. In this
figure, $y_{lm} \in \mathbb{C}$ denotes the mutual admittance between
buses $(l,m)\in \mathcal{L}$ ($y_{ml} = y_{lm}$ for all $(l,m)\in
\mathcal{L}$); $y^\text{gr}_{lm} \in \mathbb{C}$ denotes the
admittance-to-ground at end $l$ of line $(l,m)\in \mathcal{L}$; $i_{lm}
\in \mathbb{C}$ denotes current injected in line $(l,m)\in \mathcal{L}$
at bus $l$; and $\rho_{lm} \in \mathbb{C}$ denotes the ratio of the
ideal phase-shifting transformer at end $l$ of line $(l,m) \in
\mathcal{L}$ ($\rho_{lm} = 1$ if there is no transformer, the ratio is
never equal to zero). For a reference on modelling of an ideal
phase-shifting transformer, see \cite{kundur-1994}. Two ideal transformers
appear in figure \ref{fig:transmission line with pst} even though only
one or none exist per branch in a transmission network. This allows one
to describe a branch using (\ref{eq:Kirchoff's first law and Ohm's
law}).

The objective of the OPF is a second order polynomial objective function
of generated active power at each generator. Let $c_{k0},c_{k1},c_{k2}
\in \RR$ denote the coefficients of the polynomial at bus $k \in
\mathcal{G}$ as can be seen in \eqref{eq:objOPF}. These can be used to
model the cost of active generation. They can be of any value, positive
or negative, so they can also be used to model minimum deviation from a
given generation plan at each generator. Let $p_k^\text{plan}$ denote an
active generation plan at bus $k\in \mathcal{G}$. One may impose $c_{k0}
= (p_k^\text{plan})^2$ , $c_{k1} = -2p_k^\text{plan}$, and $c_{k2} = 1$
to achieve this.

\textbf{OPF:}
\begin{equation}
\min~ \sum_{k \in \mathcal{G}} c_{k2} (p^\text{gen}_k)^2 + c_{k1}p^\text{gen}_k + c_{k0}, \label{eq:objOPF}
\end{equation}
over the variables $(i_k)_{k\in \mathcal{N}},(i_{lm})_{(l,m) \in
\mathcal{L}},(p_k^\text{gen})_{k\in \mathcal{N}},(q_k^\text{gen})_{k\in
\mathcal{N}}$ and $(v_k)_{k\in \mathcal{N}}$ 
subject to
\begin{gather}
\forall\, l \in \mathcal{N},~~~ i_l = \sum_{m \in \mathcal{N}(l)} i_{lm}, \label{eq:Kirchoff's first law}
\displaybreak[0]\\
\forall\, (l,m) \in \mathcal{L},~~~   \rho_{lm}\H i_{lm} = y_{lm}^\text{gr} \frac{v_l}{\rho_{lm}} + y_{lm}(\frac{v_l}{\rho_{lm}}  - \frac{v_m}{\rho_{ml}}),
\label{eq:Kirchoff's first law and Ohm's law}
\displaybreak[0]\\
\forall\, k \in \mathcal{N} \setminus \mathcal{G},~~~ v_k i_k\H = -  p^\text{dem}_k  - \text{j} q^\text{dem}_k,
\label{eq:power demand}
\displaybreak[0]\\
\forall\, k \in \mathcal{G},~~~ v_k i_k\H = p^\text{gen}_k -  p^\text{dem}_k + \text{j}( q^\text{gen}_k - q^\text{dem}_k ),
\label{eq:power demand and generation}
\displaybreak[0]\\
\forall\, k \in \mathcal{N},~~~ p_k^{\text{min}} \leq p^\text{gen}_k \leq p_k^{\text{max}},
\label{eq:bounds on active generation}
\displaybreak[0]\\
\forall\, k \in \mathcal{N},~~~ q_k^{\text{min}} \leq q^\text{gen}_k \leq q_k^{\text{max}},
\label{eq:bounds on reactive generation}
\displaybreak[0]\\
\forall\, k \in \mathcal{N},~~~ v_k^{\text{min}} \leq |v_k| \leq v_k^{\text{max}},
\label{eq:bounds on voltage amplitude}
\displaybreak[0]\\
\forall\, (l,m) \in \mathcal{L},~~~ |v_l - v_m| \leq v_{lm}^{\text{max}},
\label{eq:bound on voltage difference}
\displaybreak[0]\\
\forall\, (l,m) \in \mathcal{L},~~~ |i_{lm}| \leq i_{lm}^{\text{max}},
\label{eq:bound on current flow}
\displaybreak[0]\\
\forall\, (l,m) \in \mathcal{L},~~~ |\text{Re}(v_l i_{lm}\H)| \leq p_{lm}^{\text{max}},
\label{eq:bound on active power flow}
\displaybreak[0]\\
\forall\, (l,m) \in \mathcal{L},~~~ |v_l i_{lm}\H| \leq s_{lm}^{\text{max}}.
\label{eq:bound on apparent power flow}
\end{gather}

Here are a few explanations for the constraints: (\ref{eq:Kirchoff's
first law}) corresponds to Kirchoff's first law; (\ref{eq:Kirchoff's
first law and Ohm's law}) corresponds to Kirchoff's first law and Ohm's
law; \eqref{eq:power demand} and \eqref{eq:power demand and generation}
correspond to power balance equations; (\ref{eq:bounds on active
generation}) corresponds to bounds on active generation; (\ref{eq:bounds
on reactive generation}) corresponds to bounds on reactive generation;
\eqref{eq:bounds on voltage amplitude} corresponds to bounds on voltage
amplitude; \eqref{eq:bound on voltage difference} corresponds to bounds
on voltage difference; \eqref{eq:bound on current flow} corresponds to
bounds on current flow; \eqref{eq:bound on active power flow}
corresponds to bounds on active power flow; and (\ref{eq:bound on
apparent power flow}) corresponds to bounds on apparent power flow.

Since the ratios of the transformers are considered fixed,
\eqref{eq:Kirchoff's first law and Ohm's law} implies that current
injected at one end of a line is a linear function of the voltages at
both ends of the line. Together with \eqref{eq:Kirchoff's first law},
this implies that there exists a complex matrix $Y$ such that $\textbf{i} = Y \vb$. This so called \textit{admittance matrix}
is defined by
\begin{equation}
Y_{lm} =
\left\{\begin{array}{cl}
\sum_{k \in \mathcal{N} \setminus \{l\}} \frac{y_{lk} + y_{lk}^\text{gr}}{|\rho_{lk}|^2} & \text{if} ~~ l = m, \\
- \frac{y_{lm}}{\rho_{ml} \rho_{lm}\H} & \text{if} ~~ (l,m)\in\mathcal{L},\\
0 & \text{otherwise.}
\end{array}\right.
\end{equation} 

\subsection{Polynomial optimization formulation of the OPF}
\label{subsec:Formulation of the OPF as a polynomial optimization problem}

In order to obtain a polynomial formulation of the OPF, we proceed in 3
steps. First, we write a formulation in complex numbers. Second, we use
it to write a formulation in real numbers. Third, we use the real
formulation to write a polynomial formulation.

\subsubsection{Formulation of the OPF in complex numbers}

Let $\textbf{a}\H$ and $A\H$ denote the conjugate transpose of a complex
vector $\textbf{a}$ and of a complex matrix $A$ respectively. It can be
deduced from~\cite{lavaei-sojoudi-2012} that there exist finite sets
$\mathcal{I}$ and $\mathcal{J}$, Hermitian matrices
$(A_k)_{k\in\mathcal{G}}$ of size $n$, complex matrices
$(B_i)_{i\in\mathcal{I}}$ and $(C_i)_{i\in\mathcal{J}}$ of size $n$, and
complex numbers $(b_i)_{i\in\mathcal{I}}$ and $(c_i)_{i\in\mathcal{J}}$
such that the OPF can be written as 
\begin{equation}
\min_{\vb\in \mathbb{C}^n} ~ \sum_{k \in \mathcal{G}} c_{k2} (\vb\H A_k \vb)^2 + c_{k1}\vb\H A_k \vb+ c_{k0},
\label{eq:complexobj}
\end{equation}
subject to
\begin{gather}
\forall\, i \in \mathcal{I},~~~ \vb\H B_i \vb~ \leq ~ b_i,
\label{eq:complexcon1}
\\
\forall\, i \in \mathcal{J},~~~ |\vb\H C_i \vb| ~ \leq ~ c_i.
\label{eq:complexcon2}
\end{gather}
Constraints
\eqref{eq:complexcon2} correspond to bounds on apparent power flow
\eqref{eq:bound on apparent power flow}. Constraints
\eqref{eq:complexcon1} correspond to all other constraints. 

\subsubsection{Formulation of the OPF in real numbers}
\label{subsubsec:Formulation of the OPF in real numbers}

Let $\xb \in \RR^{2n}$ denote $[\text{Re}(\vb)\T ~ \text{Im}(\vb)\T ]\T$
as is done in \cite{lavaei-low-2012}. In order to transform the complex
formulation of the OPF \eqref{eq:complexobj}-\eqref{eq:complexcon2} into
a real number formulation, observe that $\vb\H M \vb=(\xb\T M^\text{re}
\xb)+\text{j}(\xb\T M^\text{im}\xb)$, where the superscript $\T$ denotes
transposition,
\begin{align*}
M^\text{re}&:=
\begin{bmatrix}
\text{Re}(M) &            -\text{Im}(M) \\
\text{Im}(M) & \hphantom{-}\text{Re}(M)
\end{bmatrix},
\quad\text{and}
\\
M^\text{im}&:=
\begin{bmatrix}
\hphantom{-}\text{Im}(M) & \text{Re}(M) \\
           -\text{Re}(M) & \text{Im}(M)
\end{bmatrix}.
\end{align*}
Then \eqref{eq:complexobj}-\eqref{eq:complexcon2} becomes

\begin{equation}
\label{eq:real objective}
\min_{\xb \in \RR^{2n}} ~ \sum_{k \in \mathcal{G}} c_{k2} (\xb\T A_k^\text{re} \xb)^2 + c_{k1}\xb\T A_k^\text{re} \xb + c_{k0},
\end{equation}
subject to
\begin{gather}
\forall\, i \in \mathcal{I},~~~ \xb\T B_i^\text{re} \xb ~ \leq ~ \text{Re}(b_i),
\label{eq:real constraint1}
\\
\forall\, i \in \mathcal{I},~~~ \xb\T B_i^\text{im} \xb ~ \leq ~ \text{Im}(b_i),
\label{eq:real constraint2}
\\
\forall\, i \in \mathcal{J},~~~ (\xb\T C_i^\text{re} \xb)^2 + (\xb\T C_i^\text{im} \xb)^2 ~ \leq ~ c_i^2.
\label{eq:real constraint3}
\end{gather}

\subsubsection{Formulation of the OPF as polynomial optimization problem}
\label{subsubsec:Polynomial formulation of the OPF}

We recall that a polynomial is a function $p: \xb\in\RRn\mapsto
\sum_{\alpha\in \calA}p_\alpha \xb^\alpha$, where $\calA\subset\NNn$ is a finite
set of integer multi-indices, the coefficients $p_\alpha$ are real
numbers, and $\xb^\alpha$ is the monomial $x_1^{\alpha_1}\cdots
x_n^{\alpha_n}$. Its degree, denoted $\deg p$, is the largest
$|\alpha|=\sum_{i=1}^n\alpha_i$ associated with a nonzero~$p_\alpha$.

The formulation of the OPF in real numbers \eqref{eq:real
objective}-\eqref{eq:real constraint3} is said to be a polynomial
optimization problem since the functions that define it are
polynomials.
Indeed, the objective \eqref{eq:real objective} is a polynomial of
$\xb\in\RRnn$ of degree 4, the constraints \eqref{eq:real
constraint1}-(\ref{eq:real constraint2}) are polynomials of $\xb$ of
degree~2, and the constraints \eqref{eq:real constraint3} are
polynomials of $\xb$ of degree~4.

Formulation \eqref{eq:real
objective}-\eqref{eq:real constraint3} will however not be used below because it has
infinitely many global solutions. Indeed, formulation
\eqref{eq:complexobj}-\eqref{eq:complexcon2} from which it derives is
invariant under the change of variables $\vb\to\vb e^{\text{j}\theta}$
where $\theta \in \RR$. This invariance property transfers to \eqref{eq:real
objective}-\eqref{eq:real constraint3}. An optimization
problem with non isolated solutions is generally more difficult to solve
than one with a unique
solution~\cite{bonnans-gilbert-lemarechal-sagastizabal-2006}. This
feature manifests itself in some properties of the moment-sos approach
described in section~\ref{sec:moment-sos-approach}. For this reason, we
choose to arbitrarily set the voltage phase at bus $n$ to zero. Bearing
in mind that $v_n^\text{min} \geq 0$, this can be done by replacing
voltage constraint \eqref{eq:volt con} at bus~$n$ by \eqref{eq:phase
zero}:
\begin{gather}
(v_n^\text{min})^2 \leq x_{n}^2 + x_{2n}^2 \leq (v_n^\text{max})^2,
\label{eq:volt con} \\
x_{2n} = 0 ~~\text{and}~~ v_n^\text{min} \leq x_{n} \leq v_n^\text{max}.
\label{eq:phase zero}
\end{gather}

In light of \eqref{eq:phase zero}, a polynomial optimization problem
where there are $2n-1$ variables instead of $2n$ variables can be
formulated. More precisely, the OPF \eqref{eq:objOPF}-(\ref{eq:bound on
apparent power flow}) can be cast as the following polynomial
optimization problem

\textbf{PolyOPF:}
\begin{equation}
\min_{\textbf{x} \in \mathbb{R}^{2n-1}} ~ f_0(\xb) := \sum_{\alpha} f_{0,\alpha} \xb^\alpha,
\label{eq:obj}
\end{equation}
subject to
\begin{equation}
\forall\, i=1, \ldots, m , ~~~ f_i(\xb) := \sum_{\alpha} f_{i,\alpha} \xb^\alpha \geq 0,
\label{eq:con}
\end{equation}
where $m$ is an integer, $f_{i,\alpha}$ denotes the real coefficients of
the polynomial functions $f_i$, and summations take place
over~$\mathbb{N}^{2n-1}$. The summations are nevertheless finite because
only a finite number of coefficients are nonzero.

\section{Moment-sos approach}
\label{sec:moment-sos-approach}

We first review some theoretical aspects of the moment-sos approach (a
nice short account can be found in~\cite{anjos-lasserre-2012b}, and more
in~\cite{lasserre-2010, blekherman-parrilo-thomas-2013}). Next, we
present a set of relaxations of PolyOPF obtained by this method and
illustrate it on a simple example. Finally, we emphasize the
relationship between the moment-sos approach and the rank relaxation of
\cite{lavaei-low-2012}.

\subsection{Foundation of the moment approach}
\label{subsec:Theoretical aspects}

The moment-sos approach has been designed to find global solutions to
polynomial optimization problems. It is grounded on deep
results from real algebraic geometry. The term \textit{moment-sos}
derives from the fact that the approach has two dual aspects: the moment
and the sum of squares approaches. 
Both approaches are dual of one another in the
sense of Lagrangian duality~\cite{rockafellar-1974b}. Below, we focus on
the moment approach because it leads to SDP problems that have a close
link with the previously studied SDP method
in~\cite{lavaei-low-2012}\addnote{This approach is also the one that
provides the global minimizers.}.

\refnote{Below, we have tried to keep the presentation as simple as
possible, while preserving correctness. For example, some authors assume
that $\Kb$ must be a Borel set and that $\mu$ must be a Borel measure on
$\Kb$ (other assumptions are possible). Such terms may distract the
unfamiliar reader from the main themes of the approach, which is the
reason why we avoid them. If this approach may look rather abstract at
first glance, it highlights the origin of the SDP system
\eqref{eq:objd}-\eqref{eq:con3d} below, which is the one that is solved
in the moment-sos relaxation of order~$d$. It is also the standard way
of introducing the moment-sos approach; see for example
\cite{lasserre-2001, henrion-lasserre-loefberg-2009, lasserre-2010,
anjos-lasserre-2012b}.}Let $\Kb$ be a subset of $\mathbb{R}^{2n-1}$. The
moment approach rests on the surprising (though easy to prove) fact that
the problem $\min \{f(\xb)$: $\xb\in \Kb\}$ is equivalent to the {\em
convex\/} optimization problem
\begin{equation}
\label{lasserre-measure-pbl}
\min_{\substack{\mu ~\text{positive measure on \textbf{K}} \\\int\dmu=1}}\;\int f_0\,\dmu.
\end{equation}
Although the latter problem has a simple structure, it cannot be solved
directly, since its unknown $\mu$ is an infinite dimensional object.
Nevertheless, the realized transformation suggests that the initial
difficult global optimization problem can be structurally simplified by
judiciously expressing it on a space of larger dimension. 
The moment-sos approach goes along this way by introducing a hierarchy
of more and more accurate approximations of problem
\eqref{lasserre-measure-pbl}, hence \eqref{eq:obj}-\eqref{eq:con},
defined on spaces of larger and larger dimension.

When~$f_0$ is a polynomial and $\Kb:=\{x\in\mathbb{R}^{2n-1}$: $f_i(\xb)\geq0$, for
$i=1,\ldots, m\}$ is defined by polynomials $f_i$ like in PolyOPF, it
becomes natural to approximate the measure~$\mu$ by a finite number of
its moments.
The {\em moment\/} of $\mu$,
associated with 
$\alpha \in \mathbb{N}^{2n-1}$, is the real number $y_\alpha:=
\int \xb^\alpha\,\dmu$. Then, when~$f_0$ is the polynomial in \eqref{eq:obj},
the objective of \eqref{lasserre-measure-pbl} becomes 
$\int f_0\,\dmu =\int (\sum_{\alpha}f_{0,\alpha} \xb^\alpha)\,\dmu
=\sum_{\alpha}f_{0,\alpha} \int \xb^\alpha\,\dmu
=\sum_{\alpha}f_{0,\alpha} y_\alpha$, 
whose linearity in the new
unknown $y$ is transparent.
The constraint $\int\dmu=1$ is also readily transformed into $y_0=1$. 
In contrast, expressing which are the vectors $y$ that are moments of a
positive measure $\mu$ on~$\Kb$ (the other constraint in \eqref{lasserre-measure-pbl}) is a much more difficult task; this is
known as the {\em moment problem\/} and it is still not completely
understood in the multivariate case, despite more than a century of
work~\cite{putinar-schmudgen-2008}. It is that constraint that is
approximated in the moment-sos approach.

\subsection{Hierarchy of semidefinite relaxations}
\label{subsec:Hierarchy of semidefinite relaxations}

Lasserre \cite{lasserre-2001} proposes a sequence of relaxations for any
polynomial optimization problem like PolyOPF that grow better in
accuracy and bigger in size when the order~$d$ of the relaxation
increases. Here and below, $d$ is an integer larger than or equal to
each $v_i:=\lceil(\deg f_i)/2\rceil$ for all $i=0,\ldots,m$ (we
have denote by $\lceil\cdot\rceil$ the ceiling operator).

Let $Z \mat>= 0$ denote that $Z$ is a symmetric positive semidefinite
matrix. Define $\NNp_q:=\{ \alpha \in \NNp: |\alpha| \leq q \}$, whose
cardinality is $|\NNp_q|=\tbinom{p+q}{q}:=(p+q)!/(p!\,q!)$, and denote
by $(z_{\alpha,\beta})_{\alpha,\beta \in \NNp_q}$ a matrix indexed by
the elements of~$\NNp_q$.

\textbf{Relaxation of order d:}
\begin{equation}
\min_{(y_\alpha)_{\alpha \in \NN_{2d}^{2n-1}}} ~
\sum_{\alpha} f_{0,\alpha} y_{\alpha},
\label{eq:objd}
\end{equation}
subject to
\begin{gather}
y_0 = 1,
\label{eq:con1d}
\displaybreak[0]\\
( y_{\alpha+\beta} )_{\alpha,\beta \in \NN_{d}^{2n-1}} \mat>= 0,
\label{eq:con2d}
\displaybreak[0]\\
\forall\, i=1, \ldots, m , ~~ \sum_{\gamma} f_{i,\gamma} \left(
y_{\alpha+\beta+\gamma} \right)_{\alpha,\beta \in
\NN^{2n-1}_{d-v_i}} \mat>= 0.
\label{eq:con3d}
\end{gather}
We have already discussed the origin of \eqref{eq:objd}-\eqref{eq:con1d}
in the above SDP problem, while \eqref{eq:con2d}-\eqref{eq:con3d} are
necessary conditions to ensure that $y$ is formed of moments of some
positive measure on $\Kb$\refnote{See the claim after (1.24) in
\cite{anjos-lasserre-2012b}, which can be easily proved with the
concepts presented in the section 1.6.1 of that paper.}. When $d$
increases, these problems form a \textit{hierarchy of semidefinite
relaxations}, called that way because the objective \eqref{eq:objd} is
not affected and the feasible set is reduced. These properties show that
the optimal value of problem \eqref{eq:objd}-\eqref{eq:con3d} increases
with $d$ and remains bounded by the optimal value of
\eqref{eq:obj}-\eqref{eq:con}.

For the method to give better results, a ball constraint $\|\xb\|^2\leq M$ must be added
according to the technical assumption~1.1 in \cite{anjos-lasserre-2012b}.
For the OPF problem, this can be done easily by setting~$M$
to $\sum_{k \in \mathcal{N}} (v_k^\text{max})^2$ using \eqref{eq:bounds
on voltage amplitude} and \eqref{eq:phase zero}, without modifying the
problem. The following two properties hold in this case~\cite[theorem
1.12]{anjos-lasserre-2012b}:
 \begin{enumerate}
 \itemsep=-0.5ex
 \item
 the optimal values of the hierarchy of semidefinite relaxations
 increasingly converge toward the optimal value of
 PolyOPF\refnote{Theorem 1.12 (a) in \cite{anjos-lasserre-2012b}.},
 \item
 let $\yb^d$ denote a global solution to the relaxation of order~$d$
 and $(\eb^i)_{1\leq i\leq 2n-1}$ denotes the canonical basis of~$\NN^{2n-1}$; if PolyOPF has a unique global
 solution, then $(\yb^d_{\eb^i})_{1\leq i \leq 2n-1}$
 converges towards the global solution to PolyOPF as $d$ tends to
 $+\infty$\refnote{Theorem 1.12 (b) in \cite{anjos-lasserre-2012b}.}.
 \end{enumerate}

The largest matrix size of the moment relaxation
appears in \eqref{eq:con1d} and has the value $|\NN_{d}^{2n-1}| =
\tbinom{2n-1+d}{d}$,
where $n$ is the number of buses. For a fixed $d$, matrix size is
therefore equal to $O(n^d)$\refnote{Precisely $\tbinom{2n-1+d}{d}=
(2n-1+d)\cdots(2n)/d!=O(n^d)$.}. This makes high order relaxations too large
to compute with currently available SDP software packages. Consequently,
the success of the moment-sos approach relies wholly upon its ability to
find a global solution with a low order relaxation, for which there is
no guarantee. Note that the global solution is found by a finite order
relaxation under conditions that include the convexity of the
problem~\cite{lasserre-2008b} (not the case of PolyOPF though) or the
positive definiteness of the Hessian of the Lagrangian at the saddle
points of the Lagrangian~\cite{deklerk-laurent-2011b} (open question in
the case of PolyOPF).

\subsection{Example on a 2-bus network}
\label{s:2-bus-example}

Consider the general OPF problem presented in section
\ref{subsec:Classical formulation of the OPF} on a 2-bus network. We
will focus only on one constraint among many and write its contribution
to the first couple of relaxations of the hierarchy
described in section \ref{subsec:Hierarchy of semidefinite relaxations}.

For clarity of presentation, assume there is no apparent power flow
constraint and the objective in \eqref{eq:objOPF} is a linear function
of active power. As was remarked in section \ref{subsubsec:Polynomial
formulation of the OPF}, the degree of the objective and the degree of
the constraints of PolyOPF are thus equal to 2.
The hierarchy of semidefinite relaxations is hence defined for all
orders $d\geq1$.

Notice that since there are $n=2$ buses, the vector variable
in PolyOPF can be written as $\xb = [ x_1 ~ x_2 ~ x_3 ]$.
For clarity of presentation, assume that $v_2^\text{min} = 0$. Thus, one of the constraints of \eqref{eq:phase zero} is $x_2 \geq
0$. Based on \eqref{eq:con3d}, the expressions of
this constraint in the first and second order relaxations of the
hierarchy are \eqref{eq:order 1} and \eqref{eq:order 2} respectively:
\begin{gather}
y_{010} \geqslant 0, \label{eq:order 1} \\
\begin{bmatrix}
y_{010} & y_{110} & y_{020} & y_{011} \\
y_{110} & y_{210} & y_{120} & y_{111} \\
y_{020} & y_{120} & y_{030} & y_{021} \\
y_{011} & y_{111} & y_{021} & y_{012}
\end{bmatrix}
 \succcurlyeq 0. \label{eq:order 2}
\end{gather}
For higher orders, the size of
the matrix corresponding to the constraint grows: 10, 20, 35, etc.
Nevertheless, it is the matrix in \eqref{eq:con2d} that determines the
size of the relaxation of order $d$ as its size is greater than matrix
size in \eqref{eq:con3d}.

According to section \ref{subsec:Hierarchy of semidefinite relaxations}, vector $[ y_{100} ~ y_{010} ~ y_{001} ]$ appears in all the relaxations of the hierarchy. When  optimality is reached in the relaxations, this vector converges towards the
global solution $[ x_1^\text{opt} ~ x_2^\text{opt} ~ x_3^\text{opt} ]$ to
PolyOPF, provided it is unique (Theorem 1.12 in \cite{anjos-lasserre-2012b}). Notice that in \eqref{eq:order 2},
terms appear that correspond to monomials that do not exist in PolyOPF.
Typically, $y_{012}$ corresponds to the monomial $x_2x_3^2$ of degree 3
which is not in PolyOPF because we have restricted the degree of the
polynomials to be equal to 2.

\subsection{Moment-sos relaxations and rank relaxation}
\label{subsec:Link between hierarchy of semidefinite relaxations and rank relaxation}

When the polynomials $f_i$ defining PolyOPF are quadratic, the first
order ($d=1$) relaxation \eqref{eq:objd}-\eqref{eq:con3d} is equivalent
to Shor's relaxation~\cite{lasserre-2001b}. To make the link with the
rank relaxation of~\cite{lavaei-low-2012}, consider now the case when
the $f_i$'s are quadratic and \textit{homogeneous} like
in~\cite{lavaei-low-2012}, that is $f_i(\xb)=\xb\T A_i\xb$ for
all~$i=0,\ldots,m$, with symmetric matrices~$A_i$. Then introducing the
vector $\sb$ and the matrix $Y$ defined by $\sb_i=y_{\eb^i}$ and
$Y_{kl}=y_{\eb^k+\eb^l}$,
and $\tr$ the trace operator, the first order relaxation
reads\refnote{Indeed, when $d=1$, the constraints
\eqref{eq:con1d}-\eqref{eq:con2d} are equivalent to the positive
semidefiniteness condition in \eqref{qp-order1relation-con}. Now, since
$$
\textstyle
f_i(\xb)= \sum_{k,l}(A_i)_{kl}x_kx_l=
\sum_{k,l}(A_i)_{kl}\xb^{\eb^k+\eb^l}.
$$
the objective \eqref{eq:objd} becomes
$\sum_{k,l}(A_0)_{kl}y_{\eb^k+\eb^l}= \sum_{k,l}(A_0)_{kl}Y_{kl}=
\tr(A_0Y)$ or \eqref{qp-order1relation-obj}. As for the constraints in
\eqref{eq:con3d}, they read $0\leq \sum_{|\gamma|=2}
f_{i,\gamma}y_\gamma= \sum_{k,l}(A_i)_{kl}Y_{kl}= \tr(A_iY)$ or the
second part of \eqref{qp-order1relation-con}.}
\begin{equation}
\label{qp-order1relation-obj}
\textstyle
\min_{(\sb,Y)}\;\tr(A_0Y),
\end{equation}
subject to
\begin{equation}
\label{qp-order1relation-con}
\begin{bmatrix}
1   & \sb\T \\
\sb & Y
\end{bmatrix}
\mat>=0
\quad\mbox{and}\quad
\tr(A_iY)\geq0~~(\forall\,i=1,\ldots,m).
\end{equation}
Using Schur's complement, the positive semidefiniteness condition
in \eqref{qp-order1relation-con} is equivalent to $Y-\sb\sb\T\mat>=0$.
Since $\sb$ does not intervene elsewhere in
\eqref{qp-order1relation-obj}-\eqref{qp-order1relation-con}, it can be
eliminated and the constraints of the problem can be replaced by
\begin{equation}
\label{qp-order1relation-con'}
Y\mat>=0
\quad\mbox{and}\quad
\tr(A_iY)\geq0~~(\forall\,i=1,\ldots,m).
\end{equation}
The pair made of \eqref{qp-order1relation-obj} and
\eqref{qp-order1relation-con'} is the rank relaxation of
\cite{lavaei-low-2012}.

Here is an example of application to the OPF of the above observation:
the first order moment relaxation is equivalent to the rank relaxation
of~\cite{lavaei-low-2012} if the following conditions hold
\begin{enumerate}
 \itemsep=-0.5ex
\item the objective of the OPF \eqref{eq:objOPF} is an affine function of active power,
\item there are no constraints on apparent power flow,
\item \eqref{eq:volt con} is not replaced by \eqref{eq:phase zero} to
keep the constraints quadratic.
\end{enumerate}

\section{Numerical results}
\label{sec:Numerical results}

We present numerical results for the moment-sos approach applied to
instances of the OPF for which the rank relaxation method of
\cite{lavaei-low-2012} fails to find the global solution. We focus on
the WB2 2-bus system, LMBM3 3-bus system, and the WB5 5-bus system that
are described in \cite{bukhsh-grothey-mckinnon-trodden-2013}. Note that
LMBM3 is also found in \cite{borden-demarco-lesieutre-molzahn-2011}. For
each of the three systems, the authors of
\cite{bukhsh-grothey-mckinnon-trodden-2013} modify a bound in the data
and specify a range for which the rank relaxation fails. We consider 10
values uniformly distributed in the range in order to verify that the
rank relaxation fails and to assess the moment-sos approach.
We proceed in accordance with the discussion of section
\ref{subsec:Hierarchy of semidefinite relaxations} by adding the
redundant ball constraint. Surprisingly, the second order relaxation
whose greatest matrix size is equal to $(2n+1)n$ nearly always finds the
global solution. 

The materials used are:
\begin{itemize}
 \itemsep=-0.5ex
\item
Data of WB2, LMBM3, WB5 systems available online
\cite{bukhsh-grothey-mckinnon-trodden-2013b},
\item
Intel(R) Xeon(TM) MP CPU 2.70 GHz 7.00 Go RAM,
\item
MATLAB version 7.7 2008b,
\item
MATLAB-package MATPOWER version 3.2
\cite{murillosanchez-thomas-zimmerman-2011},
\item
SeDuMi 1.02 \cite{sturm-1999} with tolerance parameter \texttt{pars.eps}
set to $10^{-12}$ for all computations\addnote{In some cases, the
default precision $10^{-8}$ is not sufficient to obtain the global
solution at a given order~$d$, requiring to run GloptiPoly with a larger
order to get it.},
\item
MATLAB-based toolbox ``YALMIP" \cite{lofberg-2004} to compute Optimization 4
(Dual OPF) in \cite{lavaei-low-2012} that yields the solution to the
rank relaxation,
\item
MATLAB-package GloptiPoly version 3.6.1
\cite{henrion-lasserre-loefberg-2009} to compute solutions to a
hierarchy of SDP relaxations \eqref{eq:objd}-\eqref{eq:con3d}.
\end{itemize}

The same precision is used as in the solutions of the test archives
\cite{bukhsh-grothey-mckinnon-trodden-2013b}. In other words, results
are precise up to $10^{-2}$ p.u. for voltage phase, $10^{-2}$ degree for
angles, $10^{-2}$ MW for active power, $10^{-2}$ MVA for reactive power,
and cent per hour for costs. Computation time is several seconds. 

GloptiPoly can guarantee that it has found a global solution to a
polynomial optimization problem, up to a given precision. This is
certainly the case when it finds a feasible point $\xb$ giving to the
objective a value sufficiently close to the optimal value of the
relaxation.

\subsection{2-bus network: WB2}
\label{subsec:2-bus network: WB2}

\begin{figure}[H]
  \centering
    \includegraphics[width=.4\textwidth]{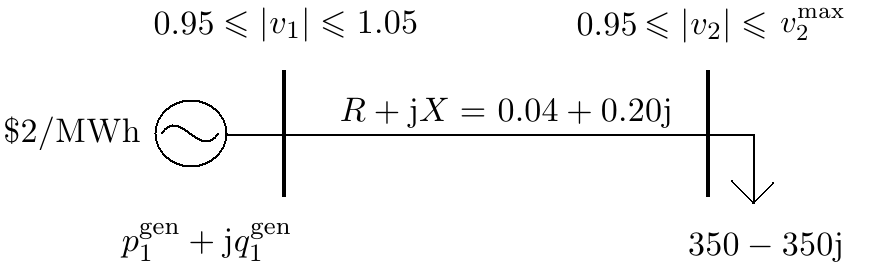}
  \caption{WB2 2-bus system}
  \label{fig:WB2 2-bus system}
\end{figure}

Authors of \cite{bukhsh-grothey-mckinnon-trodden-2013} observe that in
the WB2 2-bus system of figure \ref{fig:WB2 2-bus system}, the rank
constraint is not satisfied in the rank relaxation method of
\cite{lavaei-low-2012} when $0.976~\text{p.u.} < v_2^\text{max} <
1.035~\text{p.u.}$ In table \ref{tab:WB2}, the first column is made up
of 10 points in that range that are uniformly distributed. The second
column contains the lowest order of the relaxations that yield a global
solution. The optimal value of the relaxation of that order is written
in the third column. The fourth column contains the optimal value of the
rank relaxation (it is put between parentheses when the relaxation is
inexact).

\begin{table}[H]
\centering
\caption{Order of hierarchy needed to reach global solution to WB2 when rank relaxation fails}
\begin{tabular}{c|c|c|c}
$v_2^\text{max}$ & relax. & optimal & rank relax. \\
(p.u.) & order & value (\$/h) &  value (\$/h)\\
\hline
0.976 & 2 & 905.76 &  905.76  \\
0.983 & 2 & 905.73 & (903.12) \\
0.989 & 2 & 905.73 & (900.84) \\
0.996 & 2 & 905.73 & (898.17) \\
1.002 & 2 & 905.73 & (895.86) \\
1.009 & 2 & 905.73 & (893.16) \\
1.015 & 2 & 905.73 & (890.82) \\
1.022 & 3 & 905.73 & (888.08) \\
1.028 & 3 & 905.73 & (885.71) \\
1.035 & 2 & 882.97 &  882.97  \\
\end{tabular}
\label{tab:WB2}
\end{table}

The hierarchy of SDP relaxations is defined for $d\geq 1$ because the
objective is an affine function and there are no apparent flow
constraints. Let's explain how it works in the case where
$v_2^\text{max} = 1.022$~p.u. The optimal value of the first order
relaxation is 861.51~\$/h, that of the second order relaxation is
901.38~\$/h, and that of the third is 905.73~\$/h. This is coherent with
point~1 of the discussion of section \ref{subsec:Hierarchy of
semidefinite relaxations} that claims that the optimal values increase
with~$d$. Computing higher orders is not necessary because GloptiPoly
numerically proves global optimality for the third order.

Notice that for $v_2^\text{max} = 1.022$~p.u. the value of the rank
relaxation found in table \ref{tab:WB2} (888.08~\$/h) is different from
the value of the first order relaxation (861.51~\$/h). If we run
GloptiPoly with \eqref{eq:volt con} instead of \eqref{eq:phase zero},
the optimal value of the first order relaxation is equal 888.08~\$/h as
expected according to section \ref{subsec:Link between hierarchy of
semidefinite relaxations and rank relaxation}.

For $v_2^\text{max} = 0.976~\text{p.u.}$ and $v_2^\text{max} =
1.035~\text{p.u.}$ (see
the first and last rows of table~\ref{tab:WB2}), the rank constraint is satisfied in the rank
relaxation method so its optimal value is equal to the one of the successful
moment-sos method. In between those
values, the rank constraint is not satisfied since the optimal value is
less than the optimal value of the OPF. Notice the correlation
between the results of table \ref{tab:WB2} and the upper half of figure
8 in \cite{bukhsh-grothey-mckinnon-trodden-2013}. Indeed, the figure
shows the optimal value of the OPF is constant whereas the optimal value
of the rank relaxation decreases in a linear fashion when
$0.976~\text{p.u.} < v_2^\text{max} < 1.035~\text{p.u.}$

Surprisingly and encouragingly, according to the second column of table
\ref{tab:WB2}, the second order moment-sos relaxation finds the global
solution in 8 out of 10 times, and the third order relaxation always
find the global solution.

\textit{Remark:} The fact that the rank constraint is not satisfied for
the WB2 2-bus system of \cite{bukhsh-grothey-mckinnon-trodden-2013}
seems in contradiction with the results of papers
\cite{bose-chandy-gayme-low-2011, tse-zhang-2011, lavaei-sojoudi-2012b}.
Indeed, the authors of the papers state that the rank is less than or
equal to 1 if the graph of the network is acyclic and if load
over-satisfaction is allowed. However, load over-satisfaction is not
allowed in this network. For example, for $v_2^\text{max} = 1.022$~p.u.,
adding 1 MW of load induces the optimal value to go down from
905.73~\$/h to 890.19~\$/h. One of the sufficient conditions in
\cite{bose-chandy-gayme-low-2012} for the rank is less than or equal to
1 relies on the existence of a strictly feasible point. It is not the
case here because equality constraints must be enforced in the power
balance equation.

\subsection{3-bus network: LMBM3}

We observe that in the LMBM3 3-bus system of figure \ref{fig:LMBM3 3-bus
system}, the rank constraint is not satisfied in the rank relaxation
method of \cite{lavaei-low-2012} when $28.35~\text{MVA} \leq
s_{23}^\text{max} = s_{32}^\text{max} < 53.60~\text{MVA}$. Below
$28.35~\text{MVA}$, no solutions can be found by the OPF solver
\texttt{runopf} in MATPOWER nor by the hierarchy of SDP relaxations. At
$53.60~\text{MVA}$, the rank constraint is satisfied in the rank
relaxation method so its optimal value is equal to the optimal value of
the OPF found by the second order relaxation; see to the last row of
table \ref{tab:LMBM3}.

\begin{figure}[H]
  \centering
    \includegraphics[width=.3\textwidth]{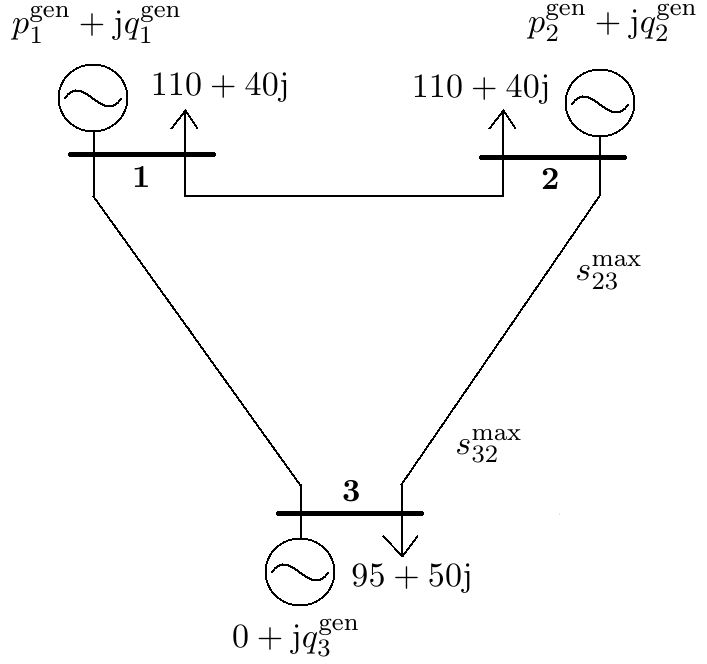}
  \caption{LMBM3 3-bus system}
  \label{fig:LMBM3 3-bus system}
\end{figure}

\begin{table}[H]
\centering
\caption{Order of hierarchy needed to reach global solution to LMBM3 when rank relaxation fails}
\begin{tabular}{c|c|c|c}
$s_{23}^\text{max} = s_{32}^\text{max}$ &  relax. & optimal & rank relax. \\
(MVA) & order & value (\$/h) & value (\$/h) \\
\hline
28.35 & 2 &            10294.88 & (6307.97) \\
31.16 & 2 & \hphantom{1}8179.99 & (6206.78) \\
33.96 & 2 & \hphantom{1}7414.94 & (6119.71) \\
36.77 & 2 & \hphantom{1}6895.19 & (6045.33) \\
39.57 & 2 & \hphantom{1}6516.17 & (5979.38) \\
42.38 & 2 & \hphantom{1}6233.31 & (5919.12) \\
45.18 & 2 & \hphantom{1}6027.07 & (5866.68) \\
47.99 & 2 & \hphantom{1}5882.67 & (5819.02) \\
50.79 & 2 & \hphantom{1}5792.02 & (5779.34) \\  
53.60 & 2 & \hphantom{1}5745.04 &  5745.04  \\ 
\end{tabular}
\label{tab:LMBM3}
\end{table}

The objective of the OPF is a quadratic function of active power so the hierarchy of SDP relaxations is defined for $d~\geqslant~2$. Again, it is surprising that the second order moment-sos relaxation always finds
the global solution to the LMBM3 system, as can be seen in the second
column of table \ref{tab:LMBM3}.

Authors of \cite{lavaei-low-2012} make the assumption that the
objective of the OPF is an increasing function of generated active
power. The moment-sos approach does not require such an assumption. For
example, when $s_{23}^\text{max} = s_{32}^\text{max} = 50~\text{MVA}$,
active generation at bus~1 is equal to 148.07~MW and active generation
at bus~2 is equal to 170.01~MW using the increasing cost function of
\cite{borden-demarco-lesieutre-molzahn-2011,
bukhsh-grothey-mckinnon-trodden-2013b}. Suppose we choose a different
objective which aims at reducing deviation from a given active generation
plan at each generator. Say that this plan is $p_1^\text{plan} =
170~\text{MW}$ at bus~1 and $p_2^\text{plan} = 150~\text{MW}$ at bus 2.
The objective function is equal to $(p_1^\text{gen}-p_1^\text{plan})^2 +
(p_2^\text{gen}-p_2^\text{plan})^2$. It is not an increasing function of
$p_1^\text{gen}$ and $p_2^\text{gen}$. The second order relaxation
yields a global solution in which active generation at bus~1 is equal to
169.21~MW and active generation at bus 2 is equal to 149.19~MW.

\subsection{5-bus network: WB5}

\begin{figure}[H]
  \centering
    \includegraphics[width=.4\textwidth]{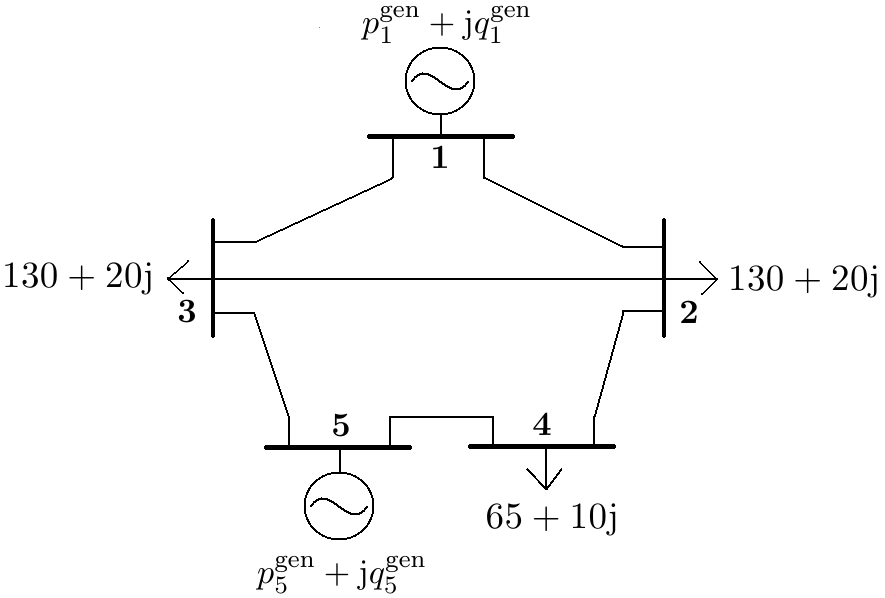}
  \caption{WB5 5-bus system}
  \label{fig:WB5 5-bus system}
\end{figure}

Authors of \cite{bukhsh-grothey-mckinnon-trodden-2013} observe that in
the WB5 5-bus system of figure \ref{fig:WB5 5-bus system}, the rank
constraint is not satisfied in the rank relaxation method of
\cite{lavaei-low-2012} when $q_{5}^\text{min}
> -30.80~\text{MVAR}$. Above $61.81~\text{MVAR}$, no solutions can be
found by the OPF solver \texttt{runopf} in MATPOWER. At
$-30.80~\text{MVAR}$, the rank constraint is satisfied in the rank
relaxation method so its optimal value is equal to the optimal value
of the OPF found by the second order moment-sos relaxation; see the
first row of table \ref{tab:WB5}. As for the 9 values considered greater
than $-30.80~\text{MVAR}$, the rank constraint is not satisfied since
the optimal value is not equal to the optimal value of the OPF. Notice that the objective of the OPF is a linear function of active power and there are bounds on apparent flow so the hierarchy of SDP relaxations is defined for $d\geqslant 1$.

\begin{table}[H]
\centering
\caption{Order of hierarchy needed to reach global solution to WB5 when rank relaxation fails}
\begin{tabular}{c|c|c|c}
$q_{5}^\text{min}$ & relax. & optimal & rank relax. \\
(MVA) & order & value (\$/h) & value (\$/h) \\
\hline
           -30.80 & 2 & \hphantom{1}945.83 & \hphantom{1(}945.83\hphantom{)} \\
           -20.51 & 2 &            1146.48 & \hphantom{1}(954.82) \\
           -10.22 & 2 &            1209.11 & \hphantom{1}(963.83) \\
\hphantom{-}00.07 & 2 &            1267.79 & \hphantom{1}(972.85) \\
\hphantom{-}10.36 & 2 &            1323.86 & \hphantom{1}(981.89) \\
\hphantom{-}20.65 & 2 &            1377.97 & \hphantom{1}(990.95) \\
\hphantom{-}30.94 & 2 &            1430.54 &            (1005.13) \\
\hphantom{-}41.23 & 2 &            1481.81 &            (1033.07) \\
\hphantom{-}51.52 & 2 &            1531.97 &            (1070.39) \\
\hphantom{-}61.81 & - &              -     &            (1114.90) \\
\end{tabular}
\label{tab:WB5}
\end{table}

When $q_{5}^\text{min} = 61.81~\text{MVAR}$, the hierarchy of SDP
relaxations is unable to find a feasible point, hence the empty slots in
the last row of table \ref{tab:WB5}. Apart from that value, the second
order moment-sos relaxation again always finds the global solution
according to the second column of \ref{tab:WB5}.

\section{Conclusion}
\label{sec:Conclusion}

This paper examines the application of the moment-sos (sum of squares)
approach to the global optimization of the optimal power flow (OPF)
problem. The result of this paper is that the OPF can be successfully
convexified in the case of several small networks where a previously known
SDP method fails. The SDP problems considered in this paper can be
viewed as extensions of the previously used rank relaxation. It is
guaranteed to be more accurate than the previous one but requires more
runtime. Directions for future research include using sparsity
techniques to reduce computational effort and identifying the OPF
problems for which a low order relaxation is exact.

\section*{Acknowledgment}

Many thanks to Stéphane Fliscounakis for fruitful discussions on
the optimal power flow problem and to Javad Lavaei for sharing his MATLAB codes with us.

 {\small
 \bibliography{article}
 \bibliographystyle{plain_e}
 }

\end{document}